%% file: REDofSYMCON.tex
\magnification=1202

\input amssym.def
\input comm

\input textstylepapernew.tex
\font\twlveusm=eusm10 at9.98pt
\font\seveneusm=eusm7 at6.38pt 
\font\fiveeusm=eusm5
\newfam\eusmfam
\textfont\eusmfam=\twlveusm
\scriptfont\eusmfam=\seveneusm
\scriptscriptfont\eusmfam=\fiveeusm
\def\eusm#1{{\fam\eusmfam\relax#1}}
\font\twlveurm=eurm10 at9.98pt
\font\seveneurm=eurm7 at6.38pt 
\font\fiveeurm=eurm5
\newfam\eurmfam
\textfont\eurmfam=\twlveurm
\scriptfont\eurmfam=\seveneurm
\scriptscriptfont\eurmfam=\fiveeurm
\def\eurm#1{{\fam\eurmfam\relax#1}}

\topskip1truecm
\voffset=2.5truecm
\hsize 15truecm
\vsize 20truecm
\hoffset=0.5truecm
\def\undertext#1{$\underline{\hbox{#1}}$}
\def\zero{\mathaccent"0017}
\topglue 3truecm
\nopagenumbers
\def\variableone{p. baguis, m. cahen}
\def\variabletwo{marsden-weinstein reduction of symplectic connections}
\def\absize{11cm}

\def\abstract#1{\baselineskip=12pt plus .2pt
                \parshape=1 0.7cm \absize
                 {\tenpbf Abstract. \tenprm #1}}
\Ref\abmar{
\author{R. Abraham, J. E. Marsden} 
\titlosb{Foundations of Mechanics}
\ekdoths{Addison-Wesley Publishing Company, Inc., pp. 299, 302 (1978)}}
				 
\Ref\baca{
\author{P. Baguis, M. Cahen}
\titlosa{A construction of symplectic connections through reduction}
\selides{\tenpit Lett. Math. Phys. {\tenprm (under press)}}}

\Ref\gotu{
\author{M. J. Gotay, G. M. Tuynman}
\titlosa{{\tenpbf R}$^{\scriptscriptstyle 2n}$
is a universal symplectic manifold for reduction}
\ekdoths{Lett. Math. Phys.}
\volume{18}
\selides{55--59 (1989)}}

\Ref\raw{
\author{J. Rawnsley}
\selides{Private communication}}

\Ref\vaisman{
\author{Vaisman, I.}
\titlosa{Connections under Symplectic Reduction}
\selides{\tenpit math.SG/0006023 v1, 4 Jun 2000}}

\hskip10cm

\centerline{\labf Marsden-Weinstein reduction}
\vskip0.2cm
\centerline{\labf for}
\vskip0.2cm
\centerline{\labf symplectic connections}

\vskip1cm

\centerline{\bf P. Baguis\footnote{$^{1}$}{e-mail: 
pbaguis@ulb.ac.be}$^{,}$\footnote{$^{2}$}{
Research 
supported by the Marie Curie Fellowship Nr. 
HPMF-CT-1999-00062}, M. Cahen\footnote{$^{3}$}{e-mail: 
mcahen@ulb.ac.be}$^{,}$\footnote{$^{4}$}{Research 
supported by an ARC of the ``Communaut\'e fran\c caise de Belgique''}}

\vskip0.3cm

\centerline{Universit\'e Libre de Bruxelles}
\centerline{Campus Plaine, CP 218 Bd du Triomphe}
\centerline{1050, Brussels, Belgium}

\vskip1cm

\abstract{We propose a reduction procedure for 
symplectic connections with symmetry. 
This is applied to coadjoint orbits whose isotropy is reductive.}

\vskip1cm

{\tenprm 
{\tenpit Key-words}: Marsden-Weinstein reduction, symplectic 
connections, Hamiltonian group actions

\vskip0.2cm

{\tenpit MSC 2000}: 53C15, 53D20}

\vfill\eject

\baselineskip=14pt plus .2pt

{\bf 0.}
The aim of this paper is to show that under very mild conditions, 
Marsden-Weinstein reduction is ``compatible'' with a symplectic 
connection. 
This means that if a symplectic manifold $(M,\omega)$ 
is endowed with a strongly 
Hamiltonian action of a connected Lie group $G$ and with a 
$G$-invariant symplectic connection $\nabla$, there is a natural way 
to construct a symplectic connection $\nabla^{r}$ on a reduced 
manifold $(M^{r},\omega^{r})$. The construction always works when $G$ 
is compact, and in many non-compact cases as well.

The interest of the construction if two-fold. First it leads to 
interesting examples of symplectic connections when $(M,\omega)$ is a 
very simple symplectic manifold and $G$ is, for example, 
one-dimensional or multidimensional but abelian (see {\baca}). 
Secondly, it may be a useful tool in dealing with the general problem 
of commutation of quantization and reduction in the framework of 
deformation quantization.

\heads

The paper is organized as follows. We first recall some classical 
results about strongly Hamiltonian actions. In the second paragraph we 
show how to construct a reduced connection with a technical assumption 
and we prove that this is always possible in the compact case. The 
third paragraph collects several examples where this construction 
gives interesting results. We finally indicate some possible further 
developments.

\vskip0.5cm

{\bf 1.} Let $(M,\omega)$ be a symplectic manifold and let 
$\sigma\colon G\times M\rightarrow M$ be a strongly Hamiltonian 
action of a connected Lie group $G$, $(g,x)\mapsto g\cdot x$,
which we will assume to be 
effective. If $\frak g$ is the Lie algebra of $G$,  we denote 
by $J\colon M\rightarrow\frak g^{\ast}$ the corresponding 
$G$-equivariant momentum map:
$$i(X^{\ast})\omega=d(J^{\ast}X),\;\forall X\an\frak g\eqn\mommap$$
where $X^{\ast}$ is the infinitesimal generator of the action 
corresponding to $X$:
$$X^{\ast}_{x}={d\over dt}\exp(-tX)\cdot x\Big|_{t=0}\eqn\infgen$$
and $J^{\ast}\colon\frak g\subset
C^{\infty}(\frak g^{\ast})\rightarrow C^{\infty}(M)$ the map defined 
by
$$(J^{\ast}X)(x)=\langle J(x),X\rangle,\;\forall x\an M.\eqn\jast$$

Let $\mu\an\frak g^{\ast}$ be a regular value of $J$ and let 
$\Sigma_{\mu}=J^{-1}(\mu)$ be the constraint manifold; it is a closed 
embedded submanifold of $M$.

The following two lemmas are classical {\abmar} and presented for 
sake of completeness.

\math{Lemma.}{\locfree}{\sl In the neighborhood of $\Sigma_{\mu}$, 
the action of $G$ is locally free, i.e. for any $x\an\Sigma_{\mu}$, 
there exists a neighborhood $\eurm\Omega_{x}$ of the identity 
element $e$ of $G$ and a neighborhood $\eurm U_{x}$ of $x$ in $M$ such 
that for any $g\an\eurm\Omega_{x}$, $y\an\eurm U_{x}$, the equation 
$g\cdot y=y$ implies $g=e$.}

\undertext{\it Proof.} Let $x\an\Sigma_{\mu}$. The map $J_{\ast 
x}\colon T_{x}M\rightarrow T_{\mu}\frak g^{\ast}\cong\frak g^{\ast}$ 
is surjective; hence the map $(J_{\ast x})^{\ast}\colon(\frak 
g^{\ast})^{\ast}\cong\frak g\rightarrow T^{\ast}_{x}M$ is injective, 
i.e. $\forall X\an\frak g$, $X\neq 0$, one has:
$$(J_{\ast x})^{\ast}(X)=(dJ^{\ast}X)_{x}=i(X^{\ast}_{x})\omega_{x}\neq 0;$$
hence $X^{\ast}_{x}\neq 0$. This means that the stabilizer $G_{x}$ 
of $x$ is discrete. Let $\chi\colon G\times M\rightarrow M\times 
M$ be the map
$(g,y)\mapsto(g\cdot y,y)$. By the above $\chi_{\ast (e,x)}$ is 
injective; hence there exist neighborhouds $\eurm\Omega_{x}$ of $e$ 
in $G$ and $\eurm U_{x}$ of $x$ in $M$ such that 
$\chi|_{\eurm\Omega_{x}\times\eurm U_{x}}$ be injective.\qed

\vskip0.3cm

Let $G_{\mu}$ be the stabilizer of $\mu$ under the coadjoint action.

\math{Lemma.}{\orthsubm}{\sl 

\item{(i)} Let $x\an\Sigma_{\mu}$ and denote by $\eusm O_{x}$ the orbit of 
$x$ under the action of $G$. Then 
$(T_{x}\Sigma_{\mu})^{\perp}=T_{x}\eusm O_{x}$ (where $\perp$ means 
orthogonal with respect to $\omega_{x}$).

\item{(ii)} Let $\Delta_{x}=(T_{x}\Sigma_{\mu})^{\perp}\cap T_{x}\Sigma_{\mu}$; 
then $\Delta_{x}$ has constant dimension (independent of $x$) and the 
orbit of $x$ under the action of $G_{\mu}$  is an integral manifold of $\Delta$.

}

\undertext{\it Proof.}

(i) For $Z\an T_{x}M$, we have:
$$J_{\ast x}Z=0\Leftrightarrow Z\an 
T_{x}\Sigma_{\mu}\Leftrightarrow\omega_{x}(X^{\ast}_{x},Z)=\langle 
J_{\ast x}Z,X\rangle=0,\;\forall X\an\frak g.$$
Consequently $T_{x}\Sigma_{\mu}\subset(T_{x}\eusm O_{x})^{\perp}$. 
But $\dim\Sigma_{\mu}=\dim M-\dim G=\dim\eusm O_{x}$ (by Lemma 
{\locfree}). Hence $T_{x}\Sigma_{\mu}=(T_{x}\eusm O_{x})^{\perp}$.

(ii) If $Z\an T_{x}\Sigma_{\mu}$, $\omega_{x}(Z,X^{\ast})=-\langle 
J_{\ast x}Z,X\rangle=0,\forall X\an\frak g$; if 
$Z\an(T_{x}\Sigma_{\mu})^{\perp}$, there exists $Y\an\frak g$ such 
that $Z=Y^{\ast}$; then $Y\an\frak g_{\mu}$, by equivariance of $J$,
where $\frak g_{\mu}$ is the Lie algebra of 
$G_{\mu}$. Hence, 
$\dim\Delta_{x}\leq\dim\frak g_{\mu}$. But (i) implies that 
$\dim\frak g_{\mu}\leq\dim\Delta_{x}$. From the above $\Delta_{x}$ is 
both the radical of $\omega|_{T_{x}\Sigma_{\mu}\times 
T_{x}\Sigma_{\mu}}$ and the tangent space to the orbit of $G_{\mu}$ 
passing through $x$. \qed

\vskip0.3cm

{\bf Assumption 1.} {\sl The constraint manifold $\Sigma_{\mu}$ is a 
$G_{\mu}$-principal bundle over the reduced manifold 
$M^{r}=G_{\mu}\backslash\Sigma_{\mu}$.}

\vskip0.2cm

{\bf Remark.} If the action of $G$ on $M$ is free and proper, 
Assumption 1 is satisfied; in particular this is true if the action is 
free and the group $G$ is compact.

\vskip0.3cm

The restriction to the constraint submanifold $\Sigma_{\mu}$ of the 
tangent bundle $TM$, denoted $TM|_{\Sigma_{\mu}}$ is a vector bundle 
over $\Sigma_{\mu}$; the group $G_{\mu}$ acts by automorphisms on this bundle. 
It contains four $G_{\mu}$-stable vector subbundles, 
$T\Sigma_{\mu}$, $(T\Sigma_{\mu})^{\perp}$, 
$T\Sigma_{\mu}+(T\Sigma_{\mu})^{\perp}$ and
$T\Sigma_{\mu}\cap(T\Sigma_{\mu})^{\perp}$.

\vskip0.3cm

{\bf Assumption 2.}{\sl There exists a $G_{\mu}$-stable vector 
subbundle $\tilde S$ of $TM|_{\Sigma_{\mu}}$ such that:
$$TM|_{\Sigma_{\mu}}=(T\Sigma_{\mu}+(T\Sigma_{\mu})^{\perp})\oplus\tilde 
S.$$}

\vskip0.2cm

{\bf Remark.} If the group $G$ is compact, such  a vector subbundle 
always exists. Indeed, we can build a $G_{\mu}$-invariant metric on 
$TM|_{\Sigma_{\mu}}$ and choose $\tilde S$ to be the orthogonal 
complement, relative to this metric, of 
$T\Sigma_{\mu}+(T\Sigma_{\mu})^{\perp}$.

\math{Lemma.}{\isotropicS}{\sl One may assume that $\tilde S$ is 
isotropic (relative to $\omega$).}

\undertext{\it Proof.} By dimension argument, $\dim\tilde S=\dim
(T\Sigma_{\mu}\cap(T\Sigma_{\mu})^{\perp})$ and $\omega$ induces a 
non-singular pairing between these two $G_{\mu}$-invariant subbundles. 
Let $x\an\Sigma_{\mu}$ and let $V_{x}$ be the symplectic subspace of 
$T_{x}M$ defined by:
$$V_{x}=\tilde S_{x}\oplus\Delta_{x}.$$
There is a unique linear map $L_{x}\colon\tilde 
S_{x}\rightarrow\Delta_{x}$ such that, $\forall u,v\an\tilde S_{x}$,
$$\omega_{x}(L_{x}u,v)=\omega_{x}(u,L_{x}v),$$
$$\omega_{x}(L_{x}u,v)+\omega_{x}(u,L_{x}v)=-\omega_{x}(u,v).$$
The graph of $L_{x}$ in $V_{x}$ is an isotropic subspace $S_{x}$ of $V_{x}$ 
such that
$$V_{x}=S_{x}\oplus\Delta_{x}.$$
Let $g\an G$; then
$$\eqalign{0&=\omega_{x}(L_{x}u,v)-\omega_{x}(u,L_{x}v)=
(g^{\ast}\omega)_{x}(L_{x}u,v)-(g^{\ast}\omega)_{x}(u,L_{x}v)\cr
\hfill&=\omega_{g\cdot x}(g_{\ast}L_{x}u,g_{\ast}v)-\omega_{g\cdot 
x}(g_{\ast}u,g_{\ast}L_{x}v)\cr}$$
$$\eqalign{-\omega_{x}(u,v)&=-(g^{\ast}\omega)_{x}(u,v)=-\omega_{g\cdot 
x}(g_{\ast}u,g_{\ast}v)\cr
\hfill&=\omega_{x}(L_{x}u,v)+\omega_{x}(u,L_{x}v)=(g^{\ast}\omega)_{x}(L_{x}u,v)
+(g^{\ast}\omega)_{x}(u,L_{x}v)\cr
\hfill&=\omega_{g\cdot x}(g_{\ast}L_{x}u,g_{\ast}v)+\omega_{g\cdot 
x}(g_{\ast}u,g_{\ast}L_{x}v).\cr}$$
By unicity, $L_{g\cdot x}=g_{\ast}\comp L_{x}\comp g^{-1}_{\ast}$ and hence the 
subbundle $S$ is $G_{\mu}$-stable.\qed

\vskip0.3cm

{\bf Remark.} By dimension argument:
$$\eqalign{(S\oplus\Delta)^{\perp}&=\left((S\oplus\Delta)^{\perp}\cap T\Sigma\right)\oplus
\left((S\oplus\Delta)^{\perp}\cap T\Sigma^{\perp}\right)\cr
\hfill&\buildrel 
not\over{=\kern-2pt=}W_{1}\oplus W_{2}\cr}$$
and the two subbundles $W_{1}$ and $W_{2}$ are $G_{\mu}$-stable.

\vskip0.5cm

{\bf 2.} We consider the situation where one has a symplectic 
manifold $(M,\omega)$, a Hamiltonian action $\sigma\colon G\times 
M\rightarrow M$ of a connected Lie group $G$ and a symplectic 
connection $\nabla$ which is $G$-invariant.

\math{Lemma.}{\compactinv}{\sl If the group $G$ is compact such a 
connection always exist.}

\undertext{\it Proof.} Let $\zero\nabla$ be any symplectic connection 
and let $X,Y$ be smooth vector fields on $M$. Define:
$$(\nabla_{X}Y)_{x}=\int_{G}\left[(g\cdot\zero\nabla)_{X}Y\right]_{x}dg=\int_{G}\left(g^{-1}_{\ast}
\zero\nabla_{g_{\ast}X}g_{\ast}Y\right)(x)dg.$$
One checks that $\nabla$ is a torsion free linear connection. 
Furthermore:
$$\omega_{x}(\nabla_{X}Y,Z)+\omega_{x}(Y,\nabla_{X}Z)=$$
$$\eqalign{\;&=
\int_{G}\left[\omega_{x}\left(g^{-1}_{\ast}\zero\nabla_{g_{\ast}X}g_{\ast}Y,Z\right)+
\omega_{x}\left(Y,g^{-1}_{\ast}\zero\nabla_{g_{\ast}X}g_{\ast}Z\right)\right]dg\cr
\hfill&=\int_{G}\left[\omega_{g\cdot 
x}\left(\zero\nabla_{g_{\ast}X}g_{\ast}Y,g_{\ast}Z\right)+\omega_{g\cdot 
x}\left(g_{\ast}Y,\zero\nabla_{g_{\ast}X}g_{\ast}Z\right)\right]dg\cr
\hfill&=\int_{G}(g_{\ast}X)_{g\cdot 
x}\omega(g_{\ast}Y,g_{\ast}Z)dg=\int_{G}X_{x}\omega(Y,Z)dg\cr
\hfill&=X_{x}\omega(Y,Z),\cr}$$
if the Haar measure $dg$ is properly normalized. \qed

\vskip0.3cm

If Assumptions 1 and 
2 are satisfied, $\Sigma_{\mu}$ (the constraint manifold) is a 
$G_{\mu}$-principal bundle over the reduced manifold $M^{r}$:
$$\pi\colon\Sigma_{\mu}\rightarrow M^{r}.$$
Furthermore, at a point $x\an\Sigma_{\mu}$, the tangent space 
$T_{x}\Sigma_{\mu}$ is the direct sum of two $G_{\mu}$-invariant 
distributions:
$$T_{x}\Sigma_{\mu}=\Delta_{x}\oplus (W_{1})_{x}$$
where $\Delta_{x}=\ker\pi_{\ast x}={\rm 
rad}^{\omega}(T_{x}\Sigma_{\mu})$. The distribution $W_{1}$ will be 
called the horizontal distribution. To $W_{1}$ is canonically 
associated a connection 1-form $\alpha$ on $\Sigma_{\mu}$ (with values 
in $\frak g_{\mu}$):
$$\alpha(U)=X,$$
if $U=\delta+w_{1}$ with $\delta_{x}=(d/dt)\exp(-tX)\cdot 
x\Big|_{t=0}=X^{\ast}_{x}$. Observe that in this framework
$$T_{x}M=\Delta_{x}\oplus(W_{1})_{x}\oplus(W_{2})_{x}\oplus S_{x}.$$
Hence we have a projection operator $P_{x}\colon T_{x}M\rightarrow 
T_{x}\Sigma_{\mu}$.

\math{Definition.}{\constrconn}{\sl If $X,Y$ are smooth vector fields, 
along $\Sigma_{\mu}$, tangent at each point to $\Sigma_{\mu}$, we 
define a linear connection $\nabla$ along $\Sigma_{\mu}$, by:
$$\nabla_{X}Y=P(\zero\nabla_{X}Y).\eqn\defconstrconn$$}

\vskip-0.5cm

\math{Lemma.}{\connpropert}{\sl $\nabla$ is a torsion free linear 
connection on $\Sigma_{\mu}$. Furthermore, $G_{\mu}$ is a group of 
affine transformations of $\nabla$.}

\undertext{\it Proof.} One has for $f\an C^{\infty}(\Sigma_{\mu})$:
$$[\nabla_{X}(fY)]_{x}=P(\zero\nabla_{X}fY)_{x}=P((Xf)Y+f\zero
\nabla_{X}Y)_{x}=(X_{x}f)Y_{x}+f(x)(\nabla_{X}Y)_{x}$$
$$\nabla_{X}Y-\nabla_{Y}X-[X,Y]=P(\zero\nabla_{X}Y-\zero\nabla_{Y}X-[X,Y])=0.$$
Also, if $Z\an\frak g_{\mu}$:
$$\eqalign{(\eusm 
L_{Z^{\ast}}\nabla)_{X}Y&=[Z^{\ast},\nabla_{X}Y]-\nabla_{[Z^{\ast},X]}Y
-\nabla_{X}[Z^{\ast},Y]\cr
\hfill&=[Z^{\ast},P\zero\nabla_{X}Y]-P\zero\nabla_{[Z^{\ast},X]}Y-P\zero\nabla_{X}
[Z^{\ast},Y]\cr
\hfill&=P\left([Z^{\ast},\zero\nabla_{X}Y]-\zero\nabla_{[Z^{\ast},X]}Y-\zero\nabla_{X}
[Z^{\ast},Y]\right)\cr}$$
using the $G_{\mu}$-invariance of $P$. Hence the conclusion as 
$\zero\nabla$ is $G_{\mu}$-invariant.\qed

\math{Lemma.}{\totgeod}{\sl The orbits of $G_{\mu}$ in $\Sigma_{\mu}$ 
are totally geodesic with respect to $\nabla$ if and only if for all 
$X,Y\an\frak g_{\mu}$ and for all vector fields $Z$ on $M$, one has:
$$\omega(P\zero\nabla_{X^{\ast}}Y^{\ast},PZ)=0.$$}
\undertext{\it Proof.} The totally geodesic condition means that 
$(\nabla_{X^{\ast}}Y^{\ast})(x)$ belongs to $\Delta_{x}$ which is the 
radical of $T_{x}\Sigma_{\mu}$.\qed

\math{Definition.}{\redconn}{\sl The reduced connection $\nabla^{r}$ 
on $M^{r}$ is defined as follows. Let $X,Y$ be smooth vector fields 
on $M^{r}$; denote by $\bar X,\bar Y$ their horizontal lifts to 
$\Sigma_{\mu}$. Then:
$$\overline{(\nabla^{r}_{X}Y)}(x)=(\nabla_{\bar X}\bar 
Y)(x)-[\alpha_{x}(\nabla_{\bar X}\bar Y)]^{\ast}.\eqn\covderhor$$}

\vskip-0.5cm

\math{Proposition.}{\redonnpropert}{\sl Formula $\covderhor$ defines a 
torsion free linear connection on $M^{r}$. Furthermore, if 
$\omega^{r}$ is the 2-form on $M^{r}$ such that 
$$\omega^{r}_{\pi(x)}(X,Y)=\omega_{x}(\bar X,\bar Y),$$
then $\omega^{r}$ is symplectic and parallel relative to $\nabla^{r}$.}

\undertext{\it Proof.} Formula $\covderhor$ defines a linear 
connection on $M^{r}$. Indeed, one has, if $g\an G_{\mu}$:
$$\eqalign{\nabla_{\bar X}\bar Y\big|_{x\cdot g}-[\alpha_{x\cdot 
g}(\nabla_{\bar X}\bar Y)]^{\ast}&=\nabla_{R_{g\ast}\bar 
X}R_{g\ast}\bar Y\big|_{x\cdot g}-R_{g\ast}\left(\Ad(g)\alpha_{x\cdot 
g}(\nabla_{\bar X}\bar Y)\right)^{\ast}\cr
\hfill&=R_{g\ast}\left[(g\cdot \nabla)_{\bar X}\bar Y\big|_{x}
-\Ad(g)\Ad(g^{-1})\alpha_{x}\left((g\cdot\nabla)_{\bar X}\bar 
Y\right)^{\ast}\right]\cr
\hfill&=R_{g\ast}\left[\nabla_{\bar X}\bar Y\big|_{x}-\alpha_{x}(\nabla_{\bar 
X}\bar Y)^{\ast}\right],\cr}$$
where $R_{g}$ is the right action which corresponds to $\sigma$: 
$R_{g}(x)=\sigma(g^{-1},x)$.
Thus formula $\covderhor$ is independent of the choice of $x$ in the 
fibre over $\pi(x)$. Also:
$$\eqalign{\overline{\nabla^{r}_{X}Y-\nabla^{r}_{Y}X-[X,Y]}&=\nabla_{\bar 
X}\bar Y-\alpha_{x}(\nabla_{\bar X}\bar Y)^{\ast}-\nabla_{\bar Y}\bar X
+\alpha_{x}(\nabla_{\bar Y}\bar X)^{\ast}-\overline{[X,Y]}\cr
\hfill&=[\bar X,\bar Y]-\alpha_{x}([\bar X,\bar 
Y])^{\ast}-\overline{[X,Y]}=0\cr}$$
and $\nabla^{r}$ is torsion free.

The 2-form $\omega^{r}$ has constant rank; furthermore, 
if $\eusm S$ denotes the cyclic sum, we have:
$$\eqalign{(d\omega^{r})_{\pi(x)}(X,Y,Z)&=\eusm 
S_{X,Y,Z}\left[X_{\pi(x)}\omega^{r}(Y,Z)-\omega^{r}_{\pi(x)}([X,Y],Z)\right]\cr
\hfill&=\eusm S_{X,Y,Z}\left[\bar X_{x}\omega(\bar Y,\bar 
Z)-\omega_{x}([\bar X,\bar Y]-\alpha_{x}([\bar X,\bar Y])^{\ast},\bar 
Z)\right]\cr
\hfill&=(d\omega)_{x}(\bar X,\bar Y,\bar Z),\cr}$$
hence $\omega^{r}$ is closed. Finally:
$$\eqalign{X_{\pi(x)}\omega^{r}(Y,Z)&=\bar X_{x}\omega(\bar Y,\bar 
Z)=\omega_{x}(\zero\nabla_{\bar X},\bar Y,\bar Z)+\omega_{x}(\bar 
Y,\zero\nabla_{\bar X}\bar Z)\cr
\hfill&=\omega_{x}(P\zero\nabla_{\bar X}\bar Y,\bar 
Z)+\omega_{x}(\bar Y,P\zero\nabla_{\bar X}\bar 
Z)=\omega_{x}(\nabla_{\bar X}\bar Y,\bar Z)+\omega_{x}(\bar 
Y,\nabla_{\bar X}\bar Z)\cr
\hfill&=\omega_{x}(\overline{\nabla^{r}_{X}Y},\bar 
Z)+\omega_{x}(\bar Y,\overline{\nabla^{r}_{X}Z})\cr
\hfill&=\omega^{r}_{\pi(x)}(\nabla^{r}_{X}Y,Z)+\omega^{r}(Y,\nabla^{r}_{X}Z),\cr}$$
which proves that $\nabla^{r}$ is symplectic.\qed

\vskip0.3cm

{\bf Formula for the curvature of the reduced connection.} Let 
$X,Y,Z$ be vector fileds on $M^{r}$. Then:
$$\eqalign{\overline{R^{r}(X,Y)Z}&=\overline{\left(\nabla^{r}_{X}\nabla^{r}_{Y}-
\nabla^{r}_{Y}\nabla^{r}_{X}-\nabla^{r}_{[X,Y]}\right)Z}\cr
\hfill&=\nabla_{\bar X}(\overline{\nabla^{r}_{Y}Z})-\alpha
\left(\nabla_{\bar X}(\overline{\nabla^{r}_{Y}Z})\right)^{\ast}-
\nabla_{\bar Y}(\overline{\nabla^{r}_{X}Z})+\alpha
\left(\nabla_{\bar Y}(\overline{\nabla^{r}_{X}Z})\right)^{\ast}\cr
\hfill&\kern12pt-\nabla_{\overline{[X,Y]}}\bar Z+\alpha
\left(\nabla_{\overline{[X,Y]}}\bar Z\right)^{\ast}\cr
\hfill&=\nabla_{\bar X}\left(\nabla _{\bar Y}\bar Z-\alpha(\nabla _{\bar 
Y}\bar Z)^{\ast}\right)-\alpha\left(\nabla_{\bar X}\left(\nabla _{\bar Y}\bar Z
-\alpha(\nabla _{\bar Y}\bar Z)^{\ast}\right)\right)^{\ast}\cr
\hfill&\kern12pt-\nabla_{\bar Y}\left(\nabla _{\bar X}\bar Z-\alpha(\nabla _{\bar 
X}\bar Z)^{\ast}\right)+\alpha\left(\nabla_{\bar Y}\left(\nabla _{\bar 
X}\bar Z
-\alpha(\nabla _{\bar X}\bar Z)^{\ast}\right)\right)^{\ast}\cr
\hfill&\kern12pt-\nabla_{[\bar X,\bar Y]-\alpha([\bar X,\bar 
Y])^{\ast}}\bar Z+\alpha\left(\nabla_{[\bar X,\bar Y]-\alpha([\bar X,\bar 
Y])^{\ast}}\bar Z\right)^{\ast}\cr
\hfill&=R(\bar X,\bar Y)\bar Z-\alpha(R(\bar X,\bar Y)\bar 
Z)^{\ast}-\nabla_{\bar X}\alpha(\nabla_{\bar Y}\bar Z)^{\ast}+\alpha
\left(\nabla_{\bar X}\alpha(\nabla_{\bar Y}\bar 
Z)^{\ast}\right)^{\ast}\cr
\hfill&\kern12pt+\nabla_{\bar Y}\alpha(\nabla_{\bar X}\bar Z)^{\ast}-\alpha
\left(\nabla_{\bar Y}\alpha(\nabla_{\bar X}\bar 
Z)^{\ast}\right)^{\ast}+\nabla_{\alpha([\bar X,\bar Y])^{\ast}}\bar 
Z\cr
\hfill&\kern12pt-\alpha\left(\nabla_{\alpha([\bar X,\bar Y])^{\ast}}
\bar Z\right)^{\ast}
\cr}$$

In the special case where $\Sigma_{\mu}$ is autoparallel with respect 
to the connection $\zero\nabla$, we have 
$\nabla_{X}Y=\zero\nabla_{X}Y$ for all vector fields $X,Y$ tangent 
to $\Sigma_{\mu}$ and the vertical subbundle in $\Sigma_{\mu}$ (which 
coincides with the radical of $\omega|_{\Sigma_{\mu}}$) is preserved 
by the connection $\nabla$. Furthermore, the reduced connection 
$\nabla^{r}$ does not depend on the choice of $S$. Indeed, for 
another subbundle $\hat S$ with the same properties as $S$, we have 
another horizontal distribution $\hat W_{1}$; if $X$ is a vector field 
on $M^{r}$, $\bar X$ and $\hat X$ its horizontal lifts with respect 
to $W_{1}$ and $\hat W_{1}$, and $\hat \alpha$ the 
connection 1-form defining $\hat W_{1}$, then $\hat X=\bar 
X+\alpha(\hat X)=\bar X-\hat\alpha(\bar X)$. If $\nabla^{\hat r}$ is 
the reduced connection defined by $\covderhor$ for the connection 
$\hat\alpha$, then one easily sees that $\widehat{\nabla^{\hat 
r}_{X}Y}=\overline{\nabla^{r}_{X}Y}-\hat\alpha(\overline
{\nabla^{r}_{X}Y})=\widehat{\nabla^{r}_{X}Y}$, which simply means 
that $\nabla^{r}$ and $\nabla^{\hat r}$ coincide. The 
reduction of the symplectic connection when $\Sigma_{\mu}$ is 
autoparallel is natural and can be performed without the machinery we 
introduce here (see $\vaisman$ for more details).

\vskip0.5cm

{\bf 3.} Coadjoint orbits are standard examples of reduced symplectic 
manifolds {\abmar} {\raw}. Let $p\colon T^{\ast}G\rightarrow G$ be the 
cotangent bundle to a connected Lie group $G$;  it can be identified,
as manifold,
to the direct product $G\times\frak g^{\ast}$ by:
$$\phi\colon T^{\ast}G\rightarrow G\times\frak 
g^{\ast},\;a\mapsto(g,L_{g}^{\ast}a),\;g=p(a),$$
where $\frak g$ is the 
Lie algebra of $G$.
The left translation by $g_{1}$ of $G$, lifts to $T^{\ast}G$ and can be 
read by the above identification, as:
$$L(g_{1})\colon G\times\frak g^{\ast}\rightarrow G\times\frak 
g^{\ast},\;(g,\xi)\mapsto(g_{1}g,\xi).$$
Similarly, the right translation by $g_{1}$ reads:
$$R(g_{1})\colon G\times\frak g^{\ast}\rightarrow G\times\frak 
g^{\ast},\;(g,\xi)\mapsto(gg_{1},\Coad(g_{1}^{-1})\xi).$$
The Liouville 1-form $\theta$ on $T^{\ast}G$, reads on $G\times\frak 
g^{\ast}$:
$$\left((\phi^{-1})^{\ast}\theta\right)_{(g,\xi)}(L_{g\ast}X+\eta)
\buildrel not\over{=\kern-2pt=}\bar\theta_{(g,\xi)}(L_{g\ast}X+\eta)=\xi(X),$$
for $X\an\frak g$, $\eta\an\frak g^{\ast}$. This gives the symplectic 
form
$$\omega_{(g,\xi)}(L_{g\ast}X+\eta,L_{g\ast}X^{\prime}+\eta^{\prime})=
\langle\eta,X^{\prime}\rangle-\langle\eta^{\prime},X\rangle-\langle\xi,[X,X^{\prime}]\rangle.$$
The fundamental vector field corresponding to the left action is
$$X^{l}(g,\xi)=-R_{g\ast}X.$$
Similarly, the fundamental vector field corresponding to the right 
action is
$$X^{r}(g,\xi)=L_{g\ast}X+\xi\comp\ad(X).$$
From this one deduces the expression of the left (resp. right) momentum 
maps:
$$J^{l}(g,\xi)=\Coad(g)\xi$$
$$J^{r}(g,\xi)=\xi.$$
If $\mu\an\frak g^{\ast}$ one constructs a constraint submanifold 
$\Sigma^{l}_{\mu}$ (resp. $\Sigma^{r}_{\mu}$) corresponding to the 
left (resp. right) action:
$$\Sigma^{l}_{\mu}=\left\{(g,\Coad(g^{-1})\mu)\;|\;g\an 
G\right\}$$
$$\Sigma^{r}_{\mu}=\{(g,\mu)\;|\;g\an G\}.$$
Let us consider the constraint manifold corresponding to the right 
action:
$$T_{(g,\mu)}\Sigma^{r}_{\mu}=\{L_{g}X\;|\; X\an\frak g\}$$
$$\left(T_{(g,\mu)}\Sigma^{r}_{\mu}\right)^{\perp}=\{X^{r}(g,\xi)\;|\;X\an\frak g\}$$
$$\left(T\Sigma_{\mu}^{r}\cap(T\Sigma_{\mu}^{r})^{\perp}\right)_{(g,\mu)}
=\{Y\an\frak g\;|\;\mu\comp\ad(Y)=0\}\cong\frak g_{\mu},$$
where $\frak g_{\mu}$ is the Lie algebra of the stabilizer $G_{\mu}$ 
of $\mu$ in the coadjoint action.
$$\left(T\Sigma_{\mu}^{r}+(T\Sigma_{\mu}^{r})^{\perp}\right)_{(g,\mu)}
=\{L_{g\ast}X+\mu\comp\ad(X)\;|\;X\an\frak g\}.$$

\math{Lemma.}{\cotangred}{\sl On $(T^{\ast}G\cong G\times\frak 
g^{\ast},\omega)$ there exists a symplectic connection $\nabla$ 
invariant by the right action of $G$.}

\undertext{\it Proof.} Let $\zero\nabla$ be the linear connection on 
$G\times\frak g^{\ast}$ defined by:
$$\zero\nabla_{\tilde{X}+\eta}(\tilde 
{X}^{\prime}+\eta^{\prime})={1\over 2 }\widetilde{[X,X^{\prime}]},$$
where the $\;\;\tilde{}\;\;$ means the corresponding left invariant vector 
field. This connection is right and left invariant but not symplectic; 
indeed, one has:
$$\eqalign{(\zero\nabla_{\tilde{X}+\eta}\omega)_{(g,\xi)}(\tilde{Y}+
\zeta,\tilde{Y}^{\prime}+\zeta^{\prime})&=(\tilde{X}+\eta)[\langle\zeta,
Y^{\prime}\rangle-
\langle\zeta^{\prime},Y\rangle-\langle\zeta,[Y,Y^{\prime}]\rangle]\cr
\hfill&\kern12pt-{1\over 2}\big(-\langle\zeta^{\prime},[X,Y]\rangle-
\langle\xi,[[X,Y],Y^{\prime}]\rangle\big)\cr
\hfill&\kern12pt-{1\over 2}\big(\langle\zeta,[X,Y^{\prime}]\rangle-
\langle\xi,[Y,[X,Y^{\prime}]]\rangle\big)\cr
\hfill&=-\langle\eta,[Y,Y^{\prime}]\rangle+{1\over 
2}\langle\zeta^{\prime},[X,Y]\rangle-{1\over 2}
\langle\zeta,[X,Y^{\prime}]\rangle\cr
\hfill&\kern12pt+{1\over 2}\langle\xi,[X,[Y,Y^{\prime}]]\rangle.
}$$
This can be projected on the space of symplectic connections as 
follows. Write
$$\nabla_{U}V=\zero\nabla_{U}V+A(U)V$$
where $A(U)$ is an endomorphism such that
$$A(U)V=A(V)U\quad\hbox{(torsion free condition)}.$$
Then choose:
$$\omega(A(U)V,W)={1\over 3}[(\zero\nabla_{U}\omega)(V,W)+
(\zero\nabla_{V}\omega)(U,W)].$$
This gives a symplectic connection which is $G$-invariant.\qed

\math{Proposition.}{\reductgroup}{\sl If the group $G_{\mu}$ is 
reductive, there exists on the reduced symplectic manifold a 
symplectic connection.}

\undertext{\it Proof.} The action of $G$ on $T^{\ast}G$ is free; hence 
Assumption 1 is satisfied. The reductiveness hypothesis ensures 
Assumption 2. \qed

\vskip0.3cm

Curvature properties of these reduced connections are worth 
investigating. We recall in particular the examples given in {\baca}. 
It seems also worthwile to read the nice Gotay-Tuynman paper {\gotu} 
thinking of connections.

\vskip1cm

{\bf Acknowledgments.} We thank our friends J. H. Rawnsley and S. Gutt 
for many useful remarks.

\vfill\eject
\vskip1cm
   \ifreferenceopen \Closeout\referencewrite \referenceopenfalse \fi
   \line{\bf\hskip0pt\hfil References\hfil}\vskip\headskip
   \vskip0.3cm
   \input referenc.txa

\end

%% file: comm.tex
 at9.98pt
\font\bb=msbm10 at9.98pt
 at9.98pt 
 at5pt
\font\cyr=wncyi10 at9.98pt
\font\eightrm=cmr8
\font\eightsc=cmcsc8
\font\labf=cmbx10 at13.1pt
 at15.74pt
\font\larm=cmr10 at13.1pt
 at15.74pt
\font\sc=cmcsc10 at9.98pt
\font\tenpbf=cmbx10 at8.32pt
\font\tenpit=cmti10 at8.32pt
\font\tenprm=cmr10 at8.32pt

\def\ad{{\rm ad}}
\def\Ad{{\rm Ad}}

\def\an{\raise0.5pt\hbox{$\kern2pt\scriptstyle\in\kern2pt$}}
\def\Ann{\hbox{\rm Ann\kern1pt}}
\def\Anns{\hbox{$\scriptstyle\rm Ann\kern0.5pt$}}

\def\arkef{\advance\chapternumber by 1\sc\roman{\the\chapternumber}}
\def\Aut{\hbox{\rm Aut\kern1pt}}
\def\bell{\hskip0pt\lower1.6pt\hbox{\bel\char'012}\kern5pt}

\def\callige{\hbox{\calligl e\kern2pt}}
\def\cheridexi{\hskip0pt\lower2pt\hbox{\cheridexia}\kern5pt}
\def\cheridexia{{\bbding\char'21}}

\def\Coad{{\rm Coad}}
\def\coker{{\rm coker\kern1pt}}
\def\Colon{\colon\kern2pt}
\def\comp{\hbox{\lower5.8pt\hbox{\larm\char'027}}}

\def\corang{{\rm corang\kern1pt}}
\def\cos{\hbox{\rm cos\kern1pt}}
\def\cosh{\hbox{\rm cosh\kern1pt}}
\def\dbaraux{\hbox{\= {\kern-2pt\= {}}}}
\def\dbar#1{\raise3pt\hbox{\dbaraux}\kern-7.8pt #1}
\def\Der{\lower0.5pt\hbox{\ygoth Der}}

\def\dim{{\rm dim\kern1pt}}
\def\double{\hbox{\kern1.5pt\bb\char'156\kern-7.6pt\char'157\kern1.5pt}}

\def\enwsh#1{{\lower2.1pt\hbox{$\buildrel{\textstyle\cup}\over
{\lower.8pt\hbox{${}_{\scriptscriptstyle#1}$}}$}}}
\def\exp{{\rm exp\kern1pt}}

\def\exten{\hbox{\callig \kern-2.5pt Ext\lower2.5pt\hbox{\kern2.5pt}}}
\def\Ham{\hbox{\rm Ham\kern1pt}}
\def\im{{\rm im\kern1pt}}
\def\k{\raise0.25pt\hbox{$\ygot k$}}

\def\ker{{\rm ker\kern1pt}}

\def\Lie{\hbox{
\callig Lie\kern2pt}}
\def\mavrodexi{\hskip0pt\lower2pt\hbox{\mavrodexia}\kern5pt}
\def\mavrodexia{{\bbding\char'15}}
\def\meriki{\hbox{\cyr\char'144\kern0.3pt}}

\def\na{\raise0.5pt\hbox{$\kern2pt\scriptstyle\ni\kern2pt$}}
\def\noan{\hbox{$\an\raise0.6pt\hbox{$\kern-6.5pt\scriptstyle
          \slash\kern3pt$}$}}

\def\oplus{\;{\mathchar"2208}\;}

\def\pounds{\rlap{\lower3.5pt\hbox{\kern2.9pt\hbox{\char'26}}}
           {\script L}}
\def\pr{\hbox{\kern3pt{\calligs p}\callig r\kern2pt}}
\def\qed{\hbox{\kern0.3cm\vrule height5pt width5pt depth-0.2pt}}
\def\QED{\hbox{\kern0.3cm\vrule height6pt width6pt depth-0.2pt}}

\def\rang{{\rm rang\kern1pt}}
\def\rank{{\rm rank\kern1pt}}

\def\san{\raise0.5pt\hbox{$\kern0.7pt\scriptscriptstyle
         \in\kern0.7pt$}}

\def\scomp{\hskip-0.05truecm\hbox{\lower5pt\hbox{$\mathchar"2017$}}
           \hskip-0.05truecm}

\def\sem{\hbox{{\script S}\kern-2.5pt\callig em\kern2pt}}

\def\sin{\hbox{\rm sin\kern1pt}}
\def\sinh{\hbox{\rm sinh\kern1pt}}

\def\styl{\hbox{\bbding\char'26}}
\def\stylo{\hskip0.3truecm\hbox{\lower1.5pt\hbox{\styl}}}
\def\times{\;{\mathchar"2202}\;}

\def\tonos{\hbox{\kern-1.3pt\lower0.7pt\hbox{$\mathchar"6013$}}}
\def\tonoskef{\hbox{$\kern-1.3pt\mathchar"6013$}}

\def\wbaraux{\hbox{\= {\kern-1.4pt\= {\kern-1.4pt\= {\kern-1.4pt\=
 {\kern-1.4pt\= {\kern-1.4pt\= {\kern-1.4pt\= {\kern-1.4pt\= {}}}}}}}}}}
\def\wbar#1{\hbox{\raise3pt\hbox{\wbaraux}\kern-30.5pt #1}}

\def\wwbaraux{\hbox{\= {\kern-1.4pt\= {\kern-1.4pt\= {\kern-1.4pt\=
{\kern-1.4pt\= {\kern-1.4pt\= {\kern-1.4pt\= {\kern-1.4pt\=
{\kern-1.4pt\= {}}}}}}}}}}}
\def\wwbar#1{\hbox{\raise3pt\hbox{\wwbaraux}\kern-34pt #1}}

\catcode`\@=11
\def\eightpoint{\eightrm}
\def\footnote#1{\edef\@sf{\spacefactor\the\spacefactor}#1\@sf
     \insert\footins\bgroup\eightpoint
     \interlinepenalty100 \let\par=\endgraf
      \leftskip=0pt \rightskip=0pt
      \splittopskip=10pt plus 1pt minus 1pt \floatingpenalty=20000
      \smallskip\item{#1}\bgroup\strut\aftergroup\@foot\let\neft}
\skip\footins=12pt plus 2pt minus 4pt
\dimen\footins=30pc

\def\line{\hbox to\hsize}

\def\title#1{\line{\hss}\line{\hss#1\hss}%
\line{\hss}\hskip-0.75truecm}

\def\author#1{{\tenprm #1:}}
\def\ekdoths#1{{\tenprm #1}}

\def\selides#1{{\tenprm #1}}
\def\titlosa#1{{\tenprm #1,}}
\def\titlosb#1{{\tenpit #1\tenprm ,}}
\def\volume#1{{\tenprm Vol. \tenpbf #1\tenprm :}}

%% file: textstylepapernew.tex
%
%
%
\def\teleia{\hbox{.}}
\newif\ifPhysRev
\def\Textindent#1{\noindent\llap{#1\enspace}\ignorespaces}
\def\nonfrenchspacing{\sfcode`\.=3001 \sfcode`\!=3000 \sfcode`\?=3000
        \sfcode`\:=2000 \sfcode`\;=1500 \sfcode`\,=1251 }
\nonfrenchspacing
\newdimen\d@twidth
 {\setbox0=\hbox{s.} \global\d@twidth=\wd0 \setbox0=\hbox{s}
        \global\advance\d@twidth by -\wd0 }
\def\removehglue{\loop \unskip \ifdim\lastskip >\z@ \repeat }
\def\roll@ver#1{\removehglue \nobreak \count255 =\spacefactor \dimen@=\z@
        \ifnum\count255 =3001 \dimen@=\d@twidth \fi
        \ifnum\count255 =1251 \dimen@=\d@twidth \fi
    \iftwelv@ \kern-\dimen@ \else \kern-0.83\dimen@ \fi
   #1\spacefactor=\count255 }
\def\step@ver#1{\relax \ifmmode #1\else \ifhmode
        \roll@ver{${}#1$}\else {\setbox0=\hbox{${}#1$}}\fi\fi }
\def\attach#1{\step@ver{\strut^{\mkern 2mu #1} }}

\normalbaselineskip = 20pt plus 0.2pt minus 0.1pt
\normallineskip = 1.5pt plus 0.1pt minus 0.1pt
\normallineskiplimit = 1.5pt
\newskip\normaldisplayskip
\normaldisplayskip = 20pt plus 5pt minus 10pt
\newskip\normaldispshortskip
\normaldispshortskip = 6pt plus 5pt
\newskip\normalparskip
\normalparskip = 6pt plus 2pt minus 1pt
\newskip\skipregister
\skipregister = 5pt plus 2pt minus 1.5pt
\newif\ifsingl@    \newif\ifdoubl@
\newif\iftwelv@    \twelv@true
\def\singlespace{\singl@true\doubl@false\spaces@t}
\def\doublespace{\singl@false\doubl@true\spaces@t}
\def\normalspace{\singl@false\doubl@false\spaces@t}
\def\Tenpoint{\tenpoint\twelv@false\spaces@t}
\def\Twelvepoint{\twelvepoint\twelv@true\spaces@t}
\def\spaces@t{\relax
      \iftwelv@ \ifsingl@\subspaces@t3:4;\else\subspaces@t1:1;\fi
       \else \ifsingl@\subspaces@t3:5;\else\subspaces@t4:5;\fi \fi
      \ifdoubl@ \multiply\baselineskip by 5
         \divide\baselineskip by 4 \fi }
\def\subspaces@t#1:#2;{
      \baselineskip = \normalbaselineskip
      \multiply\baselineskip by #1 \divide\baselineskip by #2
      \lineskip = \normallineskip
      \multiply\lineskip by #1 \divide\lineskip by #2
      \lineskiplimit = \normallineskiplimit
      \multiply\lineskiplimit by #1 \divide\lineskiplimit by #2
      \parskip = \normalparskip
      \multiply\parskip by #1 \divide\parskip by #2
      \abovedisplayskip = \normaldisplayskip
      \multiply\abovedisplayskip by #1 \divide\abovedisplayskip by #2
      \belowdisplayskip = \abovedisplayskip
      \abovedisplayshortskip = \normaldispshortskip
      \multiply\abovedisplayshortskip by #1
        \divide\abovedisplayshortskip by #2
      \belowdisplayshortskip = \abovedisplayshortskip
      \advance\belowdisplayshortskip by \belowdisplayskip
      \divide\belowdisplayshortskip by 2
      \smallskipamount = \skipregister
      \multiply\smallskipamount by #1 \divide\smallskipamount by #2
      \medskipamount = \smallskipamount \multiply\medskipamount by 2
      \bigskipamount = \smallskipamount \multiply\bigskipamount by 4 }
\def\normalbaselines{ \baselineskip=\normalbaselineskip
   \lineskip=\normallineskip \lineskiplimit=\normallineskip
   \iftwelv@\else \multiply\baselineskip by 4 \divide\baselineskip by 5
     \multiply\lineskiplimit by 4 \divide\lineskiplimit by 5
     \multiply\lineskip by 4 \divide\lineskip by 5 \fi }


\def\abstract#1{\parshape=1 0.7cm \dimen10
                {\tenpbf Abstract. \tenprm #1}}

\newcount\appendixnumber     \appendixnumber=0
\newcount\chapternumber      \chapternumber=0
\newcount\equanumber         \equanumber=0
\newcount\mathnumber         \mathnumber=0
\newcount\appequanumber      \appequanumber=0
\newcount\appmathnumber      \appmathnumber=0

\let\variableone=\relax
\let\variabletwo=\relax
\let\chapterlabel=\relax
\let\sectionlabel=\relax
\let\mathlabel=\relax
\newtoks\chapterstyle        \chapterstyle={\Number}
\newtoks\sectionstyle        \sectionstyle={\chapterlabel\Number}
\newskip\chapterskip         \chapterskip=\bigskipamount
\newskip\sectionskip         \sectionskip=\medskipamount
\newskip\headskip            \headskip=8pt plus 3pt minus 3pt
\newdimen\chapterminspace    \chapterminspace=15pc
\newdimen\sectionminspace    \sectionminspace=10pc
\newdimen\sectionspace       \sectionspace=20pc
\newdimen\referenceminspace  \referenceminspace=25pc

\def\chapterreset{\global\advance\chapternumber by 1
   \ifnum\equanumber<0 \else\global\equanumber=0\fi
   \mathnumber=0
   \makechapterlabel}
\def\makechapterlabel{\let\sectionlabel=\relax\let\mathlabel=\relax
 \xdef\chapterlabel{\the\chapterstyle{\the\chapternumber\teleia\kern3pt}}}

\def\rightheadline{\sc\hfil\variableone\eightsc\hfil\folio}
\def\leftheadline{\eightsc\folio\hfil{\sc\variabletwo}\hfil}
\def\heads{\footline={\hfil}\headline={\ifodd\pageno
               \rightheadline\else\leftheadline\fi}}

\def\headseis{\partreset\headline={\ifodd\pageno{
                         \hfil\sc partie {\eightsc\the\partnumber}
                         -introduction\hfil\eightsc\folio}\else
                        {\eightsc\folio\hfil\sc partie
                         {\eightsc\the\partnumber}-introduction\hfil}\fi}
                        \footline={\hfil}}

\def\alphabetic#1{\count255='140 \advance\count255 by #1\char\count255}
\def\Alphabetic#1{\count255='100 \advance\count255 by #1\char\count255}
\def\Roman#1{\uppercase\expandafter{\romannumeral #1}}
\def\roman#1{\romannumeral #1}
\def\Number#1{\number #1}
\def\BLANC#1{}

\def\titlestyle#1{\par\begingroup \interlinepenalty=9999
     \leftskip=0.02\hsize plus 0.23\hsize minus 0.02\hsize
     \rightskip=\leftskip \parfillskip=0pt
     \hyphenpenalty=9000 \exhyphenpenalty=9000
     \tolerance=9999 \pretolerance=9000
     \spaceskip=0.333em \xspaceskip=0.5em
     \iftwelv@\bf\else\bf\fi
   \noindent #1\par\endgroup }

\def\spacecheck#1{\dimen@=\pagegoal\advance\dimen@ by -\pagetotal
   \ifdim\dimen@<#1 \ifdim\dimen@>0pt \vfil\break \fi\fi}
\def\TableOfContentEntry#1#2#3{\relax}

\def\chapter#1{\par\vskip0.7cm
   \chapterreset \titlestyle{\chapterlabel\ #1}
   \nobreak\vskip\headskip
   \wlog{\string\chapter\space \chapterlabel} }

\def\appendixreset{\global\advance\appendixnumber by 1
                   \appmathnumber=0\appequanumber=0}
\def\appendix#1{\par \penalty-300\vskip\chapterskip
   \spacecheck\chapterminspace
   \appendixreset \title{\bf Appendix \Alphabetic{\the\appendixnumber}}
   \nobreak\vskip-\chapterskip\penalty 30000
   \vskip-\chapterskip
   \par{\titlestyle{#1}}
   \vskip\chapterskip
   \wlog{\string\appendix\space \chapterlabel} }

%
%
\def\eqname#1{\relax \ifnum\equanumber<0
     \xdef#1{{\noexpand\rm(\number-\equanumber)}}%
       \global\advance\equanumber by -1
    \else \global\advance\equanumber by 1
      \xdef#1{{\noexpand(
                             \rm{\number\equanumber})}} \fi #1}

\def\eqn{\eqno\eqname}

\def\math#1#2{\vskip0.1cm
   \global\advance\mathnumber by 1
   \xdef\mathlabel{\the\mathnumber}
   \wlog{\string\math\space \mathlabel}
   {\bf\enspace\mathlabel\hskip0.2cm #1}
   \xdef#2{{\mathlabel}}}

\def\appeqname#1{\relax \ifnum\appequanumber<0
     \xdef#1{{\noexpand\rm(\number-\appequanumber)}}%
       \global\advance\appequanumber by -1
    \else \global\advance\appequanumber by 1
      \xdef#1{{\noexpand(\hbox{\Alphabetic{\the\appendixnumber}}\teleia
                            {\number\appequanumber})}} \fi #1}

\def\mathapp#1#2{\vskip0.1cm
   \global\advance\appmathnumber by 1
   \xdef\appmathlabel{{\Alphabetic{\the\appendixnumber}}\teleia
   \the\appmathnumber}
   \wlog{\string\mathapp\space \appmathlabel}
   {\bf\enspace\appmathlabel\hskip0.2cm #1}
   \xdef#2{{\appmathlabel}}}


%
%
%
\newtoks\referencestyle      \referencestyle={\tenpbf\Number}
\newcount\referencecount     \referencecount=0
\newcount\lastrefsbegincount \lastrefsbegincount=0
\newif\ifreferenceopen       \newwrite\referencewrite
\newif\ifrw@trailer
\newdimen\refindent     \refindent=13pt
\def\NPrefmark#1{\attach{\scriptscriptstyle [ #1 ] }}
\let\PRrefmark=\attach
\def\refmark#1{\relax\ifPhysRev\PRrefmark{#1}\else\NPrefmark{#1}\fi}
\def\refend@{\refmark{\number\referencecount}}
\def\refend{\refend@{}\space }
\def\refsend{\refmark{\count255=\referencecount
   \advance\count255 by-\lastrefsbegincount
   \ifcase\count255 \number\referencecount
   \or \number\lastrefsbegincount,\number\referencecount
   \else \number\lastrefsbegincount-\number\referencecount \fi}\space }
\def\refitem#1{\par\hangafter=0 \hangindent=\refindent	\Textindent{#1}}
\def\Ref{\rw@trailertrue\REF}
\def\REF#1{\r@fstart{#1}%
   \rw@begin{\tenprm [\tenpbf\Number{\the\referencecount}\tenprm ]}\rw@end}
\def\r@fstart#1{\chardef\rw@write=\referencewrite \let\rw@ending=\refend@
   \ifreferenceopen \else \global\referenceopentrue
   \immediate\openout\referencewrite=referenc.txa
   \toks0={\catcode`\^^M=10}\immediate\write\rw@write{\the\toks0} \fi
   \global\advance\referencecount by 1 
   \xdef#1{[{\the\referencestyle{\the\referencecount}}]}}
 {\catcode`\^^M=\active %
 \gdef\rw@begin#1{\immediate\write\rw@write{\noexpand\refitem{#1}}%
   \begingroup \catcode`\^^M=\active \let^^M=\relax}%
 \gdef\rw@end#1{\rw@@end #1^^M\rw@terminate \endgroup%
   \ifrw@trailer\rw@ending\global\rw@trailerfalse\fi }%
 \gdef\rw@@end#1^^M{\toks0={#1}\immediate\write\rw@write{\the\toks0}%
   \futurelet\n@xt\rw@test}%
 \gdef\rw@test{\ifx\n@xt\rw@terminate \let\n@xt=\relax%
       \else \let\n@xt=\rw@@end \fi \n@xt}%
}
\let\rw@ending=\relax
\let\rw@terminate=\relax

\def\Closeout#1{\toks0={\catcode`\^^M=5}\immediate\write#1{\the\toks0}%
   \immediate\closeout#1}
%